\documentclass[12pt,a4paper]{amsart}

\usepackage{amsmath,amssymb}
\usepackage{bbm,paralist}
\usepackage{a4wide}

\newtheorem{theorem}{Theorem}
\newtheorem{proposition}[theorem]{Proposition}
\newtheorem{lemma}[theorem]{Lemma}

\theoremstyle{definition}

\newcommand{\RR}{\mathbb{R}}

\newcommand{\ov}{\overline{v}}
\newcommand{\oV}{\overline{V}}

\newcommand{\e}{\varepsilon}
\newcommand{\be}{\boldsymbol{\e}}
\newcommand{\one}{\boldsymbol{1}}

\newcommand{\matrixrule}{\raisebox{.5ex}{\rule{1.5em}{.4pt}}}

\DeclareMathOperator{\proj}{proj}

\DeclareMathOperator{\Id}{I}

\DeclareMathOperator{\vertices}{vert}

\begin{document}

\title[Explicit neighborly centrally symmetric polytopes]{An explicit construction for neighborly\\ centrally symmetric polytopes}
\author%% [%
%% {Julian Pfeifle}
%% \qquad\MakeUppercase{\textbf{\sffamily --- \today{} ---}}
%% ]
{Julian Pfeifle}

\address{Departament de Matem\`{a}tica Aplicada II, Universitat
  Polit\`{e}cnica de Catalunya}

\email{julian.pfeifle@upc.edu}

\thanks{The author was supported by a \emph{Juan de la Cierva} postdoctoral
  research grant from the Spanish Ministry of Education (MEC), and by projects
  MCYT BFM2003-00368 and MTM2005-08618-C02-01}

\date{June 29, 2006}

\begin{abstract}
  We give an explicit construction, based on Hadamard matrices, for an
  infinite series of $\big\lfloor\frac12\sqrt{d}\big\rfloor$-neighborly
  centrally symmetric $d$-dimensional polytopes with $4d$~vertices.
  This appears to be the best explicit version yet of a recent
  probabilistic result due to Linial and Novik, who proved the existence of
  such polytopes with a neighborliness of~$\frac{d}{400}$.
\end{abstract}

\maketitle

\section{Introduction}

A polytope $P\subset\RR^d$ is \emph{centrally symmetric} (cs, for short) if
$P=-P$. A cs~polytope~$P$ is \emph{$k$-neighborly} if every set of $k$~of its
vertices, no two of which are antipodes, is the vertex set of a face of~$P$.

In their recent paper~\cite{Linial-Novik06}, Linial and Novik give
probabilistic constructions for highly neighborly cs~polytopes. Namely, based
on probabilistic techniques due to Garnaev and
Gluskin~\cite{Garnaev-Gluskin84}, they construct $k$-neighborly
$d$-dimensional cs~polytopes with $2m=2(n+d)$ vertices, such that
$k=\Theta\big(\frac{d}{1+\log(m/d)}\big)$; moreover, they show that this value
is asymptotically optimal.  In the ``diagonal'' case $n=d$ they use a
probabilistic result due to Ka\v sin~\cite{Kasin77} to construct
$d$-dimensional $\frac{d}{400}$-neighborly cs~polytopes with $4d$~vertices,
and ask if there exists an explicit construction of highly neighborly
cs~polytopes.

%\smallskip
In this note, we provide such an explicit and non-probabilistic
construction:

\begin{theorem}\label{thm:thm}
  For each $d\ge4$ such that there exists a Hadamard matrix of size~$d$, there
  is an explicit construction for a
  $\big\lfloor\frac12\sqrt{d}\big\rfloor$-neighborly cs~$d$-polytope with
  $4d$~vertices.
\end{theorem}

Hadamard matrices exist  for every $d=2^e$ with integer $e\ge2$;
see~\cite{Seberry-Yamada92} for a survey.

%\smallskip
Of course, Theorem~\ref{thm:thm} does not attain (by far) the bound given by
Linial and Novik's probabilistic arguments, but to date no better explicit
construction of highly neighborly cs~polytopes seems to be known; see also our
additional comments in Section~\ref{sec:discuss}.  We refer
to~\cite{Linial-Novik06} for a (necessarily short) survey of the known results
on neighborly cs~polytopes.

\smallskip To briefly outline the remaining contents of this note, we proceed
to characterize $k$-neighborly cs~$d$-polytopes with $2(n+d)$~vertices in
terms of a certain linear projection (Proposition~\ref{prop:proj}) and a
certain matrix equation (Theorem~\ref{thm:matrix}). In
Section~\ref{sec:hadamard}, we then find very special solutions of this
equation in the ``diagonal'' case~$n=d$ and prove Theorem~\ref{thm:thm}.

\section{cs-transforms, polarity, and a matrix equation}

The key to our construction is the following lemma due to Linial and Novik,
which characterizes those point sets arising as McMullen and Shephard's
\emph{cs-transforms}~\cite{McMullen-Shephard68} of cs~polytopes:

\begin{lemma}[Linial and Novik~\mbox{\cite[Lemma~3.1]{Linial-Novik06}}]
  \label{lem:1}
  A cs~set $\oV=\{\pm \ov_1,\dots, \pm\ov_m\}\subset\RR^n$ is a cs~transform
  of the vertex set of a $k$-neighborly cs~$d$-polytope with $2m=2(d+n)$
  vertices if and only if the set $\oV_+=\{\ov_1,\dots, \ov_m\}$ does not
  contain dominant subsets of size~$k$.
\end{lemma}
  
Here, Linial and Novik define a subset $\{\ov_i:i\in I\}$ of $\oV_+$ to be
\emph{dominant} if there exists $0\ne u\in \RR^n$ such that
\begin{equation}\label{eq:dominant}
  \sum_{i\in I} \big|\langle\ov_i,u\rangle\big| \ \ge \ 
  \frac12 \sum_{j=1}^m \big|\langle\ov_j,u\rangle\big|\,. 
\end{equation}
  
To interpret this characterization geometrically, let
$\Diamond_m\subset\RR^m$~be the standard cross-polytope, $C_{\pm a}(m)$ the
$m$-dimensional cube $\{x\in\RR^m: \mbox{$-a\le x_i\le a$} \text{ for } 1\le
i\le m\}$, and $E(m,k)$ the $m$-dimensional convex hull of all
$0/\mbox{$\pm1$}$-vectors of length~$m$ with exactly $k$~non-zero entries.
Thus, $E(m,k)=C_{\pm1}(m)\cap k\Diamond_m$; equivalently, we obtain $E(m,k)$
by reflecting the standard $m$-dimensional hypersimplex~$\Delta(m,k)$ in the
coordinate hyperplanes of~$\RR^m$.
  
Now let $T:\RR^n\to\RR^m$ be the linear map given by the real $(m\times
n)$-matrix whose rows are $\ov_1,\dots,\ov_m$, so that
$Tu=(\langle\ov_1,u\rangle, \dots, \langle\ov_m,u\rangle)^T$ for any
$u\in\RR^n$ (here and throughout, the superscript~$T$ denotes transpose; we
trust that this will not cause confusion). We assume that $T$~has full rank,
and denote the image of~$T$ by~$L$, a linear $n$-space in~$\RR^m$.
  
To express Linial \& Novik's lemma in this language, write
$\be=(\e_1,\dots,\e_m)$ and note~that
\begin{eqnarray*}
   \frac12\sum_{j=1}^m \big|\langle \ov_j,u \rangle\big| 
   & = &
   \frac12\max\left\{
     \sum_{j=1}^m \e_j \langle \ov_j,u \rangle\ : \
     \be\in\vertices C_{\pm1}(m)\right\}\\
   & = &
   \frac12\max\left\{
     \be
     \big(\langle \ov_1,u \rangle, \dots, \langle \ov_m,u \rangle\big)^T
      :  \be\in\vertices C_{\pm1}(m)   
   \right\}\\
   & = &
   \max\left\{ 
     \langle z, Tu \rangle : z\in\vertices C_{\pm\frac12}(m) 
   \right\},
\end{eqnarray*}
and analogously, for any subset $I\subset[m]$
of cardinality $k$,
\begin{eqnarray*}
   \sum_{i\in I} \big|\langle \ov_i,u \rangle\big| 
   & = &
   \max\left\{
     \sum_{i\in I}\delta_i \langle \ov_i,u \rangle \ : \
     \delta_i=\pm1 \text{ for all } i\in I
   \right\}\\
   & = &
   \max\big\{ 
      \langle w,Tu \rangle \ : \ w \in\vertices E(m,k)
     \big\}. 
\end{eqnarray*}
Thus, Lemma~\ref{lem:1}, condition~\eqref{eq:dominant} and the fact that the
maximum of any linear function on a polytope is attained at one of the
vertices  together say that $\{\pm\ov_1,\dots,\pm\ov_m\}$ is a
cs-transform of a $k$-neighborly cs~$d$-polytope with $2m$~vertices if and
only if
\[
      \max_{w\in E(m,k)} \langle w,v \rangle \ < \ 
      \max_{z\in C_{\pm1/2}(m)} \langle z,v \rangle
      \qquad\text{for all }  0\ne v\in L\,.
\]
By dualizing~---~i.e., considering $\langle v,w\rangle$ instead of $\langle
w,v\rangle$~---~we can also read this condition as saying that for
\emph{any} non-zero vector $v\in L$, an affine hyperplane perpendicular to~$v$
that sweeps outward from the origin along~$v$ should have left behind all
vertices of~$E(m,k)$ before encountering the last vertex of~$C_{\pm 1/2}(m)$.
We have reached the following conclusion:

\begin{proposition}\label{prop:proj}
  The set $\{\pm\ov_1,\dots,\pm\ov_m\}\subset\RR^n$ is a cs-transform of a
  $k$-neighborly cs~$d$-polytope with $2m$~vertices if and only if
  \begin{equation}\label{eq:proj} 
    \proj_L E(m,k) \ \subset \ \proj_L C_{\pm1/2}(m),
  \end{equation}
  where $\proj_L$ denotes orthogonal projection to~$L\subset\RR^m$, the linear
  $n$-space that is the image of the linear map $T:\RR^n\to\RR^m$ whose matrix
  has rows $\ov_1,\dots,\ov_m$.
\end{proposition}

To proceed, we take advantage of the following duality (Lemma~\ref{lem:polar})
that relates the section of a polytope~$P\subset\RR^m$ by a linear
subspace~$L$ to the projection of~$P^\Delta$ to~$L$. Recall that the
\emph{polar set} of~$P$ is
\begin{equation}\label{eq:defpolar}
      P^\Delta \ = \ \big\{x\in \RR^m: \langle x,y\rangle\le 1 \text{ for all
      } y\in P\big\},
\end{equation}
and that $(P^\Delta)^\Delta=P$ if $0\in P$.

\begin{lemma}\label{lem:polar}
  Let $P\subset\RR^m$ be any polytope such that $0\in P$, let
  $L\subset\RR^m$~be any linear subspace, and denote the orthogonal projection
  of~$\RR^m$ to~$L$ by~$\proj_L$. Then
    \begin{equation}\label{eq:polar}
      \proj_L(P^\Delta) \ = \ (P\cap L)^\Delta\cap L.
    \end{equation}
\end{lemma}
  
We learned about this lemma from~\cite{Ziegler-comm}; the proof is elementary
and follows from the definition~\eqref{eq:defpolar} of a polar set.

%%  \begin{proof}
%%    Via a suitable linear transformation of the ambient space~$\RR^m$, we may
%%    assume $L$ to be the linear subspace of codimension~$d$ given by
%%    $x_{n+1}=\dots=x_m=0$, where $n=m-d$. We must show the equality of the two
%%    sets
%%      \begin{eqnarray*}
%%         \proj_L (P^\Delta) & = & \Big\{ (x_1,\dots,x_n,0,\dots,0)\in\RR^m : 
%%         \text{ there exist } x_{n+1},\dots,x_m\in\RR \\
%%         && \phantom{\big\{} \text{ such that } 
%%         \langle x,y \rangle \le 1 \text{ for all } y\in P, \text{ where }
%%         x=(x_1,\dots,x_m)\Big\}
%%      \end{eqnarray*}
%%      and
%%      \[
%%         (P\cap L)^\Delta \cap L \ = \ 
%%         \big\{ x = (x_1,\dots,x_n,0,\dots,0)\in\RR^m: 
%%         \langle x,y\rangle\le 1 \text{ for all } y\in
%%         P\cap L\big\}.
%%      \]
%%      The forward inclusion follows because any $y\in P\cap L$ is of the form
%%      $y=(y_1,\dots,y_n,0,\dots,0)$.  
%%  \end{proof}

By  substituting~\eqref{eq:polar} with $P=E(m,k)^\Delta$, respectively
$P=C_{\pm1/2}(m)^\Delta$, into~\eqref{eq:proj}, we obtain
\[
   \big( E(m,k)^\Delta\cap L \big)^\Delta \cap L \ \subset \ 
   \big( C_{\pm 1/2}(m)^\Delta\cap L \big)^\Delta \cap L.
\]
We now restrict to the subspace~$L$ and polarize. Because both polytopes
contain the origin and are full-dimensional in~$L$ (and therefore polarizing
reverses inclusion, and $(P^\Delta)^\Delta=P$), we obtain the equivalent
condition
\[
      C_{\pm1/2}(m)^\Delta \cap L \ \subset \ E(m,k)^\Delta \cap L.
\]
This in turn is satisfied if and only if the $2^m$ facet-defining inequalities
of the polytope $C_{\pm1/2}(m)^\Delta=2\Diamond_m$, together with some fixed
set of $d=m-n$ equations defining~$L$, imply the $2^k\binom{m}{k}$
facet-defining inequalities of~$E(m,k)^\Delta$.

To find a linear subspace $L$ that achieves this, we represent $L$ as the
kernel of the matrix $(\Id_d|A)$, where $\Id_d$~is the $(d\times d)$ identity
matrix and $A=(a_{ij})$~a real $(d\times n)$ matrix. Moreover, we pass to
homogeneous coordinates, which means to express each point $x\in\RR^m$ as
$(1,x)\in\RR\times\RR^m$, and each inequality $ax\le a_0$, for $a\in(\RR^m)^*$
and $a_0\in\RR$, as $(a_0,a)\in (\RR^{m+1})^*$. 

Phrased in this language, we must express each vertex~$(1,e)$ of $\{1\}\times
E(m,k)\in \RR\times\RR^m$ as a linear combination of the following form:
\begin{equation}\label{eq:comb}
    \begin{tabular}[c]{c|*{7}{c}}
      $\mu_1^e$ & $2$ & $-\delta_{1,1}^e$  & $\cdots$ & $-\delta_{1,d}^e$ & 
      $-\delta_{1,d+1}^e$ & $\cdots$ & $-\delta_{1,d+n}^e$ \\
      $\cdots$ &&& $\cdots$ &&& $\cdots$ \\
      $\mu_{2^m}^e$ & $2$ & $-\delta_{2^m,1}^e$  & $\cdots$ &
      $-\delta_{2^m,d}^e$ &  
      $-\delta_{2^m,d+1}^e$ & $\cdots$ & $-\delta_{2^m,d+n}^e$ \\
      $\nu_1^e$ & 0 & 1 & $\cdots$ & 0 & $a_{1,1}$ & $\cdots$ & $a_{1,n}$\\
      $\cdots$ &&& $\cdots$ &&& $\cdots$ \\
      $\nu_d^e$ & 0 & 0 & $\cdots$ & 1 & $a_{d,1}$ & $\cdots$ &
      $a_{d,n}$\\
      $\e$ & 1&0&$\cdots$&0&0&$\cdots$&0\\
      \hline\hline\\[-2ex]
      & 1 &&&& \hspace{-1.5cm}$e$ 
    \end{tabular}
 \end{equation}
  
 \noindent In this table, the $\mu$'s, $\nu$'s and $\e$ are understood to
 multiply the adjacent row vectors, and the result of this linear combination
 is the row vector $(1,e)$. Specifically, $\e\ge0$, $\mu^e_i\ge0$, and
 $\nu_i^e\in\RR$ for all relevant indices, and not all of these coefficients
 are required to be non-zero; moreover, $\delta_{ij}^e=\pm1$, and the
 reason for introducing the minus signs for the~$\delta$'s will become clear
 in a moment.  We will also use the notation $M^e_j=\sum_{i\in
   I_e}\mu^e_i\delta^e_{ij}$, for $1\le j\le d+n$, where
 $I_e\subset\{1,\dots,2^m\}$ indexes the non-zero $\mu^e_i$.  Note the
 constraint $\sum_{i\in I_e}\mu^e_i\le\frac12$ implied by the ``$0$-th''
 column of this linear combination, which in turn implies $|M^e_j|\le\frac12$.
  
 From columns $1\le j\le d$ of~\eqref{eq:comb}, we learn that $\sum_{i\in
   I_e}\mu^e_i(-\delta^e_{ij})+\nu^e_j=e_j$, so that in fact we know the
 coefficients $\nu^e_j=e_j+M^e_j$.  With this information, we obtain from
 columns $d+1\le j\le d+n$ that $\sum_{i=1}^d(e_i+M^e_i)a_{i,j-d}=e_j+M^e_j$.
 Expressed in matrix notation, we have arrived at the following result:

 \begin{theorem}\label{thm:matrix}
   Finding a cs-transform of a $k$-neighborly cs $d$-polytope with
   $2m=2(n+d)$ vertices is equivalent to finding a $(d\times n)$-matrix~$A$
                                %of rank~$d$ 
   and a $(2^k\binom{m}{k} \times m)$-matrix~$M'$ that satisfy the
   following requirements:
   \begin{enumerate}[\upshape(a)]
   \item All entries of~$M'$ are bounded in absolute value by~$\frac12$.
   \item Let $E'=(E|F)$ be a $(2^k\binom{m}{k}\times m)$-matrix whose rows are
     the vertices of~$E(m,k)$ in some order, and decompose it into a
     matrix~$E$ with $d$~columns and a matrix~$F$ with $n$~columns. Similarly,
     decompose $M'=(M|N)$ into a matrix~$M$ with $d$ columns and a matrix~$N$
     with $n$~columns.
     Then the matrices $A$, $E$, $F$, $M$ and $N$ must satisfy
     \begin{equation}\label{eq:matrixeq}
       (E+M)\,A \ = \ F+N.
     \end{equation}
   \item The rows of $M'$ must be expressible as linear combinations as
     in~\eqref{eq:comb}.
   \end{enumerate}
 \end{theorem}
  
 As an aside, it is clear that any expression of~$(1,e)$ as a linear
 combination as in~\eqref{eq:comb} immediately yields an expression
 of~$(1,-e)$ as a similiar linear combination, by reversing the signs of the
 relevant $\delta^e_{ij}$ and $\nu_i$.  It would therefore be enough to
 consider only one member of each pair of antipodal vertices of~$E(m,k)$, and
 consequently only keep those rows of the $\{0,\pm1\}$-matrix~$E'$ whose first
 non-zero entry is positive.  However, to keep the symmetry of the problem we
 choose not to do this.
 
 We do, however, partially order the rows of~$E'$. Namely, we partition~$E'$
 into $k+1$ blocks $E_l'=(E_l|F_{k-l})$ with $0\le l\le k$, such
 that each row of~$E_l$ has exactly $l$~non-zero entries (and
 consequently each row of~$F_{k-l}$ has $k-l$~of them). The order inside
 each such block is immaterial for our purposes. Note that with this partial
 ordering, the number of rows of both~$E_l$~and~$F_{k-l}$ is
 \[
     2^{l}\binom{d}{l}\cdot 2^{k-l}\binom{n}{k-l} \ = \
     2^{k}\binom{d}{l}\binom{n}{k-l}.
 \]
 By decomposing $M$ and $N$ into blocks labeled $M_{k-l}$, $N_{l}$ with the
 same number of rows as $E_l$~and~$F_{k-l}$, equation~\eqref{eq:matrixeq}
 above decomposes into the $k+1$ equations
 \begin{equation}
    \label{eq:matrixeq2}
    (E_l+M_{k-l})\,A \ = \ F_{k-l} + N_{l}\,,
    \qquad \text{ for \ $0\le l\le k$. }
 \end{equation}

\section{Hadamard matrices}\label{sec:hadamard}
  
In the diagonal case $n=d$, we will exhibit a very special solution
of~\eqref{eq:matrixeq2}. Namely, we find matrices $A$,
$M_{k-l}$,~$N_{l}$ such that
 \begin{eqnarray*}
    E_l A &=& N_{l}, \\
    M_{k-l} A &=& F_{k-l},
  \end{eqnarray*}
in the following way:
  
Let $d$ be such that there exists a Hadamard matrix of order~$d$, i.e., a
$(d\times d)$ matrix~$H_d$ with entries~$\pm1$ such that $H_d^TH_d=d\Id_d$,
and set
\[
     A \ = \ \alpha H_d \ = \ 
     \alpha
     \begin{pmatrix}
       \vrule & & \vrule\\ v_1 & \cdots &v_d\\ \vrule && \vrule
     \end{pmatrix}
     \ = \ 
     \alpha
     \begin{pmatrix}
       \matrixrule & w_1 & \matrixrule \\ & \vdots\\ \matrixrule  & w_d &
       \matrixrule  
     \end{pmatrix},
 \]
 for some real constant $\alpha>0$ to be determined later. Here the
 $(\pm1)$-vector~$v_i$ denotes the $i$-th column, and respectively $w_j$ the
 $j$-th row of~$H_d$; this implies that $v_i^Tv_i=d$ and $v_i^Tv_j=0$ for
 $j\ne i$, and similarly for the $w$'s.
%% , so that
%%  $|v_{ij}|=|w_{ij}|=\alpha$ for $1\le i,j\le n$. 
 Moreover, set 
 \[
    N_{l} \ = \ E_l A
 \]
 and
 \[
     M_{k-l} \ = \ 
     \beta F_{k-l} H_d^T \ = \
     \beta
     \begin{pmatrix}
       v_1^T+\dots+v_{k-l}^T\\
       v_1^T+\dots-v_{k-l}^T\\
       \dots\\
       -v_{d-k+l+1}^T-\dots-v_d^T
     \end{pmatrix},
 \]
 for $\beta>0$ another real constant. The displayed pattern of signs and
 indices in~$M_{k-l}$ reflects the one in~$F_{k-l}$, and thus corresponds to
 a fixed but arbitrary ordering of the rows of~$F_{k-l}$.
  
 We now adjust $k$, $\alpha$ and $\beta$ to make these matrices compatible
 with the conditions in Theorem~\ref{thm:matrix}. For this, first note that
 each row of $N_{l} = E_l A$ is of the form $\alpha\sum_{i\in
   I}\sigma_iw_i$, for some index set $I\in\binom{[d]}{l}$ and signs
 $\sigma_i\in\{\pm1\}$. In particular, the absolute value of each entry
 of~$N_{l}$ is bounded by $\alpha l\le \alpha k$,
 so that $\alpha$~is constrained by
 \begin{equation}\label{eq:star}
   \alpha k \ \le \ \frac12.
 \end{equation}
 Similarly, each element of~$M_{k-l}$ is bounded in absolute value by
 $\beta(k-l)\le\beta k$, so we also need
 \begin{equation}\label{eq:starstar}
   \beta k \ \le \ \frac12. 
 \end{equation}
 Because $A=\alpha H_d$, we obtain $M_{k-l} A=\alpha\beta F_{k-l}H_d^T H_d =
 \alpha\beta d F_{k-l}$, so we must set $\beta=\frac{1}{\alpha d}$ in order to
 fulfill~\eqref{eq:matrixeq2}, and thus condition~(b) of
 Theorem~\ref{thm:matrix}. Now \eqref{eq:star}~and~\eqref{eq:starstar} taken
 together say that $k,\alpha$~and~$d$ must satisfy
 \begin{equation}\label{eq:kcond}
     \frac{2k}{d} \ \le \ \alpha \ \le \ \frac{1}{2k}, 
 \end{equation}
 so that we arrive at the bound $k\le\frac12\sqrt{d}$ for the
 cs-neighborliness of our cs-polytope.  In fact, for $d\ge4$, the choices
 $k:=\big\lfloor\frac12\sqrt{d}\big\rfloor$, $\alpha:=\frac{1}{2k}$ and
 $\beta:=\frac{1}{\alpha d}=\frac{2k}{d}$ satisfy
 \eqref{eq:star}~and~\eqref{eq:starstar}, and thus we have found a matrix
 $M'=(M|N)$ that satisfies conditions~(a)~and~(b) of Theorem~\ref{thm:matrix}.
  
 It only remains to check condition~(c), i.e., that the rows of~$M'$ can in
 fact be expressed as linear combinations as in~\eqref{eq:comb}.  For this,
 note that by the definitions of~$M_{k-l}$~and~$N_{l}$, each row of~$M'$
 is a sum of row vectors of the form
 \[
    \bigg(\alpha\sum_{i\in I} \sigma_i v_i^T \, \bigg| \, \boldsymbol0\bigg) +
    \bigg(\boldsymbol0 \, \bigg| \, \beta\sum_{j\in J} \sigma_j w_j\bigg)
 \]
 for index sets $I\in\binom{[d]}{k-l}$ and $J\in\binom{[d]}{l}$ and
 signs $\sigma_i,\sigma_j=\pm1$, where $\boldsymbol0$ denotes the zero row
 vector of length~$d$ and $0\le l\le k$.  We now represent
 \begin{eqnarray}
     \Big(\alpha \sum_{i\in I} \sigma_i v_i^T\Big|\boldsymbol0\Big) & = & 
     \sum_{i\in I}\Big(
     \tfrac{\alpha}{2} \big(\sigma_i v_i^T \big| \boldsymbol1\big)+
     \tfrac{\alpha}{2} \big(\sigma_i v_i^T \big| {-\boldsymbol1}\big)\Big),
     \label{eq:comb1}\\
     \Big(\boldsymbol0\Big|\beta \sum_{j\in J} \sigma_j w_j\Big) & = & 
     \sum_{j\in J}\Big(
     \tfrac{\beta}{2}\big(\boldsymbol1 \big| \sigma_j w_j\big)+
     \tfrac{\beta}{2}\big({-\boldsymbol1} \big| \sigma_j w_j\big)\Big)
     \label{eq:comb2}
\end{eqnarray}
as linear combinations of $2(k-l)$, respectively $2l$, vectors of length
$2d$ with entries~$\pm1$, where $\boldsymbol1$~represents the all-ones
vector of length~$d$.  The sum over all coefficients in this linear
combination is then
\begin{eqnarray*}
  \sum_{i\in I}\left(\frac{\alpha}{2}+\frac{\alpha}{2}\right) +
  \sum_{j\in J}\left(\frac{\beta}{2}+\frac{\beta}{2}\right) & = &
  (k-l)\frac{1}{2k} + l\,\frac{2 k}{d}\\
  &=& \frac12 - l\left(\frac{1}{2k}-\frac{2k}{d}\right)\\
  &\le& \frac12
\end{eqnarray*}
by~\eqref{eq:kcond}, as required. This concludes the proof of
Theorem~\ref{thm:thm}.

\section{Discussion}\label{sec:discuss}

We are plainly still quite far away from an explicit construction of
$\Theta(d)$-neighborly $d$-dimensional cs-polytopes with $4d$~vertices. This
situation is all too familiar: Linial \& Novik find the linear subspace~$L$
defined by our matrix~$A$ using a probabilistic construction due to Ka\v sin,
and remark on the difficulty of explicitly finding such subspaces. 

In the light of the discussion in Ball~\cite[p.~24]{Ball97}, our explicit
construction of $\Theta(\sqrt{d})$-neighborly cs~polytopes using Hadamard
matrices is what can reasonably be expected in this context, and it may not be
realistic to hope for more: \emph{``There are some good reasons, related to
  Ramsey theory, for believing that one cannot expect to find genuinely
  explicit matrices of any kind that would give the right estimates''}.

We close the present note by briefly mentioning some variations and
alternatives.

\subsection{Special Hadamard matrices}
The bound $k=O(\sqrt{d})$ arises via~\eqref{eq:star},~\eqref{eq:starstar}
from~\eqref{eq:kcond} because $m$~is our best a priori upper bound for the
largest absolute value of an entry of the sum of $m$~rows (or columns)
of~$H_d$. If this largest absolute value could instead be taken of
order~$O(\sqrt{m})$ for~$m=O(d)$, we would reach our goal of a
cs-neighborliness of $k=\Theta(d)$.

To address this issue, our construction of Section~\ref{sec:hadamard} works
for \emph{any} Hadamard matrix, but these are in fact quite a varied and
structured lot, cf.~\cite{Craigen-Kharaghani04}. In particular, there exist
so-called \emph{regular} Hadamard matrices of order~$d$, for which all the
entries of the sum of \emph{all}~$d$ rows (or columns) are
precisely~$\sqrt{d}$.  However, this is not good enough for our purposes: it
follows from elementary considerations that any row or column of a
(conveniently normalized) regular Hadamard matrix contains exactly
$\frac12(d+\sqrt{d})$ entries~`$1$' and $\frac12(d-\sqrt{d})$ entries~`$-1$';
therefore, there exist choices of $l=O(d)$ rows or columns such that the
maximal entry of their sum will be $O(d)$~in absolute value, and
via~\eqref{eq:kcond} this ruins our cs-neighborliness.

Another reason for doubting the efficacy of Hadamard matrices in this respect
is that the fraction of the total number of vertices of the $2d$-dimensional
cube involved in the concrete instances~\eqref{eq:comb1},~\eqref{eq:comb2} of
the linear combination~\eqref{eq:comb} is quite small.

\subsection{Pseudo-inverses} 
Moving away from Hadamard matrices, one should really try to \emph{find} the
right matrix~$A$ in~\eqref{eq:matrixeq}~or~\eqref{eq:matrixeq2}, instead of
prescribing it. In this context, we recall the concept of \emph{generalized
  inverses}, and refer to~\cite{Ben-Israel03} for further discussion and
notation. 

A \emph{Moore-Penrose $\{1\}$-inverse} of a
real $(m\times n)$~matrix~$G$ is any real $(n\times m)$~matrix~$G^{(1)}$ such
that $GG^{(1)}G=G$; the set of all $\{1\}$-inverses of~$G$ is
denoted~$G\{1\}$.  These are important for our purposes because
by~\cite[Theorem 2.1]{Ben-Israel03}, the matrix equation 
\[
   GA \ = \ H
\]
has a solution $A$ if and only if there exists  $G^{(1)}\in G\{1\}$
such that
\[
   (GG^{(1)}-I_m) H \ = \ 0.
\] 

To apply these notions to our context, we set $G=E+M$ and $H=F+N$, and remark
that in the diagonal case $n=d$, an especially nice ordering of the rows
of~$E'_l=(E_l|F_{k-l})$ is the ``doubly lexicographic'' one:

 \begin{proposition}\label{rem:orthog}
   If $n=d$, one can choose a total order on the rows of~$E'$ (that refines
   the partial order given above), in such a way that the matrices
   $E_l$~and~$F_{k-l}$ satisfy
   \begin{eqnarray*}
     E_l^T  E_l & = & 2^k\binom{d-1}{l-1}\binom{d}{k-l}\, I_d,\\
     E_l^T  F_{k-l} & = & 0.
   \end{eqnarray*}
\end{proposition}

\begin{proof}
  Let $\tilde C_l$~be the matrix of size $\binom{d}{l}\times d$ whose rows
  are, in lexicographical order, all $0/1$-vectors of length~$d$ with exactly
  $l$~entries~`$1$', set $n(l)=2^l\binom{d}{l}$, and let $C_l$~be the
  matrix of size $n(l)\times d$ obtained from $\tilde C_l$ by replacing each
  row with the $2^l$~rows obtained by choosing all possible signs for the
  non-zero entries, again in lexicographical order. Thus, the non-zero entries
  of each column of~$C_l$ come in $\binom{d-1}{l-1}$~blocks of size~$2^l$
  each, so that the scalar product of each column with itself
  is~$2^l\binom{d-1}{l-1}$. Moreover, it easily follows by induction that
  distinct columns of~$C_l$ are mutually orthogonal, so that
   \[
      C_l^T C_l \ = \ 2^l\binom{d-1}{l-1} I_d.
   \]
   Now denote the all-ones column vector of length~$i$ by~$\one_i$. Then 
   \begin{eqnarray*}
       E_l & = & C_l \otimes \one_{n(k-l)},\\
       F_{k-l} &=& \one_{n(l)} \otimes C_{k-l}
   \end{eqnarray*}
   combine in such a way that the matrix $(E_l|F_{k-l})$ is a valid
   representation of~$E'_l$. (Recall that $A\otimes B$ is the matrix obtained
   from~$A$ by replacing each entry~$a_{ij}$ by the block matrix $a_{ij}B$, so
   that $(A\otimes B)^T = A^T\otimes B^T$ and $(A\otimes B)(C\otimes D) = AC
   \otimes BD$.) Now, as claimed,
   \[
      E_l^T E_l \ = \ 
      ( C_l^T \otimes \one_{n(k-l)}^T ) ( C_l \otimes \one_{n(k-l)} ) \ = \ 
      2^k \binom{d-1}{l-1} \binom{d}{k-l} I_d,
   \]
   and
   \[
      E_l^T F_{k-l} \ = \
      ( C_l^T \otimes \one_{n(k-l)}^T ) ( \one_{n(l)} \otimes C_{k-l} )
      \ = \ ( C_l^T \one_{n(l)} ) \otimes ( \one_{n(k-l)}^T C_{k-l} )
      \ = \ 0,
   \]
   because, again by induction, the sum of all entries in any column
   of each~$C_i$  vanishes.
 \end{proof}

Therefore, we can choose our matrices $G$ and $H$ to be
\[
   G \ = \
   \begin{pmatrix}
     C_{k\phantom{{}-1}} \phantom{{}\otimes \one_{n(k-1)}} 
     + M_{k\phantom{{}-1}} \\
     C_{k-1} \otimes \one_{n(1)\phantom{{}-k}} + M_{k-1} \\
     \dots \\
     C_{1\phantom{{}-k}} \otimes \one_{n(k-l)} + M_{1\phantom{{}-k}} \\
     \phantom{C_{1} \otimes \one_{n(k-1)} +{}} M_{0}
   \end{pmatrix}, 
   \qquad
   H \ = \
   \begin{pmatrix}
     \phantom{C_{k-1} \otimes \one_{n(1)} +{}} N_{0} \\
     \one_{n(k-1)} \otimes C_{1\phantom{{}-k}}  + N_{1\phantom{{}-k}} \\
     \dots\\
     \one_{n(1)\phantom{{}-k}} \otimes C_{k-1} + N_{k-1} \\
     \phantom{\one_{n(k-1)} \otimes{}} 
     C_{k\phantom{{}-1}} + N_{k\phantom{{}-k}}
   \end{pmatrix}.
\]
The set of all $\{1\}$-inverses of~$G$ can be parametrized
explicitly using the techniques in~\cite{Udwadia-Kalaba99}; however, so far we
have not succeeded in turning this to our advantage.

\end{document}